\documentclass[12pt]{article}
\usepackage{amssymb}
\usepackage{amscd}
\usepackage{epsf}
\newtheorem{thm}{Theorem}[section]

\newtheorem{prop}[thm]{Proposition}

\newcommand{\inv}{^{-1}}

\newcommand{\Z}{{\mathbb{Z}}}
\newcommand{\R}{{\mathbb{R}}}

\newcommand{\PP}{{\mathbb{P}}}
\newcommand{\Gm}{{\mathbb{G}}_m}

\newcommand{\cs}{{\cal{S}}}

\newcommand{\cn}{{\cal {N}}}

\newcommand{\ov}{\overline}

\newcommand{\ga}{\gamma}

\newcommand{\bg}{{\bf  G}}
\newcommand{\bn}{{\bf N}}
\newcommand{\bt}{{\bf T}}

\newcommand{\n}{\underline{n}}
\newcommand{\La}{\Lambda}
\newcommand{\Labar}{\overline{\Lambda}}
\newcommand{\abar}{\overline{A}}
\newcommand{\xbar}{\overline{X}}

\newcommand{\etabari}{\overline{\eta}_i^I}

\newcommand{\an}{^{an}}

\begin{document}

\title{Compactification of the Bruhat-Tits building of PGL 
by seminorms}
{ \author{Annette Werner \\ \small Mathematisches Institut, 
Universit\"at M\"unster, Einsteinstr. 62, D -  48149 M\"unster\\
 \small e-mail: werner@math.uni-muenster.de}}

\date{}

\maketitle

\centerline{\bf Abstract}
\small 
We construct a compactification $\overline{X}$
of the Bruhat-Tits building $X$ associated 
to the group $PGL(V)$ which can be identified with the space of
homothety classes of seminorms on $V$ endowed with the topology
of pointwise convergence. Then we define a continuous map
from the projective space to $\overline{X}$ which extends
the reduction map from Drinfeld's $p$-adic symmetric domain to the 
building $X$.

\centerline{{\bf 2000 MSC:} 20E42, 20G25} 

\normalsize 
\section*{Introduction}

It is well-known that the Bruhat-Tits building $X$ of $PGL(V)$ can 
be identified with the set of homothety classes of norms on $V$.
It can also be described with lattices of full rank in $V$.
In \cite{we} we constructed a compactification of $X$ which 
takes into account lattices
of arbitrary rank in $V$.

The goal of the present paper is to describe a related compactification
$\xbar$,  which is more adapted to the description of the building in terms
of norms. It is in some sense dual to the construction given in \cite{we}.
Whereas there the boundary of $X$ consists of all buildings corresponding 
to subspaces of $V$, here we attach the buildings corresponding to 
quotient spaces of $V$. 
We prove  that
our construction leads to a compact, contractible space $\xbar$, and 
 that $PGL(V)$ acts continuously on $\xbar$. 

A nice feature of this new construction is that it has a very
natural description: 
Namely, $\xbar$ can be identified with the set
of homothety classes of seminorms, endowed with the topology of pointwise
convergence.

In the final section, we use the compactification $\xbar$ in order to 
show that the reduction map $r : \Omega\an \rightarrow X$ from the Berkovich
analytic space corresponding to Drinfeld's $p$-adic symmetric domain to the 
building $X$ has a natural extension to a continuous, $PGL(V)$-equivariant
map
\[ r: \PP(V)\an \longrightarrow \xbar,\]
where $\PP(V)\an$ is  the Berkovich analytic space induced by the projective
space. Besides we show
that $r$ has a continuous section $j: \xbar \rightarrow \PP(V)\an$, which 
induces a homeomorphism between $\xbar$ and a closed subset of $\PP(V)\an$. 

We hope that the compactification $\xbar$ can be used to construct
a kind of Satake compactification for Bruhat-Tits buildings associated
to arbitrary reductive groups.

{\bf Acknowledgements: }I am much indepted to Matthias Strauch for 
some very fruitful discussions. In fact, the idea to modify the construction 
in \cite{we} 
in order to extend the reduction map from Drinfeld's symmetric domain
to the building was born during one of them. 
Besides, it is a pleasure to thank Vladimir Berkovich for
some useful conversations about his analytic spaces.

\section{Notation and conventions}
Throughout this paper we denote by $K$ a non-archimedean local field, 
by $R$ its valuation ring and by $k$ the residue class field. Besides, 
$v$ is the valuation map, normalized so that it maps a prime element to $1$, 
$q$ is the number of elements in the residue field and $|x|=q^{-v(x)}$ the
absolute value on $K$.

We adopt the convention that ``$\subset$'' always means strict subset, 
whereas we write ``$\subseteq$'', if equality is permitted.

Let $V$ be an $n$-dimensional vector space over $K$, and let $\bg$ be the algebraic
group $ \underline{PGL}(V)$.

\section{Seminorms}

We call a map 
$\gamma: V \rightarrow \R_{\geq 0 }$
a seminorm, if $\gamma$ is not identically zero, and satisfies

i) $\gamma(\lambda v) = |\lambda| \gamma(v)$ for all $\lambda \in K$ and $v \in V$, and

ii) $\gamma(v+w) \leq \sup\{ \gamma(v), \gamma(w)\}$ for all $v,w$ in $V$.

We call $\gamma$ canonical with respect to a basis $v_1, \ldots, v_n$ of $V$, if
\[ \gamma(\lambda_1 v_1 + \ldots + \lambda_n v_n ) = \sup \{|\lambda_1| \gamma(v_1),
 \ldots, |\lambda_n| \gamma(v_n)\} \]
  for all
$\lambda_1, \ldots, \lambda_n$ in $K$.

A seminorm $\ga$ satisfying 
 
 iii) $\gamma(x) > 0 $ for all $x \neq 0$,
 
 is a norm on $V$. 
 
A seminorm $\gamma$ on $V$ induces in a natural way a norm on the quotient space 
$V/\mbox{ker} \gamma$. By \cite{goi}, Proposition 3.1, for any norm $\gamma$ 
there exists
a basis with respect to which $\gamma$ is canonical. Looking at the norms 
induced on a quotient, we find that the same holds
for seminorms. 

Two norms or seminorms $\gamma_1$ and $\gamma_2$ on $V$ are called equivalent,
iff there is a positive real constant $c$ such that $\ga_1 = c \ga_2$.

\section{Compactification of  one apartment}
We fix a basis $v_1,\ldots, v_n$ of $V$, which defines a maximal $K$-split
torus $\bt$ in $\bg$, induced by the torus $\bt^\sim$ 
of diagonal matrices in $GL(V)$ with
respect to $v_1,\dots, v_n$. We denote by $\chi_i$ the character of $\bt^\sim$
mapping a diagonal matrix to its i-th entry. Then all $\chi_i/\chi_j$ define
a character $a_{ij}$ of $\bt$. 

Let $\bn$ be the normalizer of $\bt$ in $\bg$. 
We denote by $T= \bt(K)$, $N = \bn(K)$ and $G = \bg(K)$ the groups of 
rational points. Then $W = N/T$ is the Weyl group of the root system 
$\Phi = \{ a_{ij}: i \neq j\} $
corresponding to $\bt$. We can identify $W$ in a natural way with the
group of permutations of $\{1,\ldots, n\}$. By embedding $W$ as the group
of permutation matrices in $N$, we find that $N$ is the semidirect product
of $T$ and $W$.

By $X_\ast(\bt)$ respectively $X^\ast(\bt)$ we denote the cocharacter 
respectively the character group of $\bt$. Let $A$ be the 
$\R$-vector space $A= X_\ast(\bt) \otimes_\Z \R$. We take $A$
as our fundamental appartment in the definition of the Bruhat-Tits building
associated to $\bg$. 

Let $\eta_i : \Gm  \rightarrow \bt$ be the cocharacter induced by 
mapping $x$ to the diagonal matrix with diagonal entries $d_1,\ldots,d_n$ 
such that $d_k= 1$ 
for $k \neq i$ and $d_i = x$. Then $\eta_1 + \dots + \eta_n = 0$, and
 $\eta_1,\ldots, \eta_{n-1}$ is an $\R$-basis of $A$.
If $t \in T$ is induced by the diagonal matrix with
entries $d_1, \ldots, d_n$, we define a point $\nu(t) \in A$
by $\nu(t) = -v(d_1) \eta_1 - \ldots  - v(d_n)  \eta_n$.

Every $t \in T$ acts on $A$ 
as translation by $\nu(t)$,
Besides,  $W$ acts as
a group of reflections on $A$. Since $N$ is the semidirect product
of $T$ and $W$, we can define an action of $N$ on $A$ by affine
bijections. 

It is well-known that $A$ can be identified with the 
set $\cn_c$ of equivalence classes of those norms on $V$
which are canonical with respect to the basis $v_1,\ldots, v_n$, cf. \cite{brti}.
To be more explicit, define a map 
\[\varphi: A \longrightarrow \cn_c\]
by mapping 
 $x = x_1 \eta_1 + \ldots + x_n \eta_n$ 
in $A$ to the class of norms represented by 
\[\ga(\lambda_1 v_1 + \ldots + \lambda_n v_n) = \sup \{ |\lambda_1| q^{-x_1} , 
\ldots, |\lambda_n| q^{-x_n}\}.\]
It is easily seen that $\varphi$ is well-defined, bijective and  $N$-equivariant, 
if we let $N$ act on 
the class of equivalence classes of norms by 
$\gamma \mapsto \gamma \circ n\inv$.

 Let us denote
by $\cs_c'$ the set of all seminorms on $V$ which are canonical with
respect to $v_1,\ldots, v_n$, and by $\cs_c$ the quotient of $\cs_c'$ by
the equivalence relation on seminorms defined above.  
We will now define a compactification $\abar$ of $A$ so that 
$\varphi$ can be extended to a homeomorphism from $\abar$ to $\cs_c$. 

We write $\n$ for the set $\{1, \ldots, n\}$.
Let $J$ be a non-empty subset of $\n$, and let $V_J$ be the subspace of $V$
 generated by the $v_i$ for $i \in J$. We write $\bg^{V_J}$ for the subgroup of 
 $\bg = \underline{PGL}(V)$ consisting of the elements leaving the subspace $V_J$
 invariant, and 
 $\bg_J$ for the algebraic group $\underline{PGL}(V_J)$. Then we have a natural restriction map
$ \rho_J: \bg^{V_J} \longrightarrow \bg_J$.
The torus $\bt$ is contained in  $\bg^{V_J}$, and its image under $\rho_J$ 
is a maximal $K$-split torus $\bt_J$ in $\bg_J$, namely the torus induced 
by the  diagonal matrices with respect to the base $\{v_i : i \in J\}$ of $V_J$.

Besides, we have  the quotient map 
$q_J: V \rightarrow V /V_J$. There is a natural homomorphism
\[\sigma_{\n \backslash J} : \bg^{V_J} \rightarrow
\underline{PGL}(V/V_J) .\]

Then for all subsets $I \subset \n$ the homomorphism
$\sigma_I$ maps $\bt$ to a split torus $\ov{\bt}_I$ in 
$\underline{PGL}(V/V_{\n \backslash
I})$. Put $
\overline{T}_I =  \overline{\bt}_{I}(K)$. Besides, 
we put
 \[A_{I} = X_\ast(\overline{\bt}_{I}) \otimes \R .\]
Then $\sigma_I$ induces 
a surjective homomorphism 
\[s_{I}: A \rightarrow A_{I},\] 
mapping the cocharacter
$\eta_i$ to zero, if $i$ is not contained in $I$, and to the induced cocharacter
$\etabari$ of $\overline{\bt}_I$, if $i $ is contained in $I$. Then $A_{I}$ 
is generated
by the cocharacters $\etabari$ for $i \in I$,
subject to the relation $\sum_{i \in I} \etabari = 0$.

Besides,  we define a homomorphism
$\ov{\nu}_I:  \ov{T}_I \rightarrow A_{I}$ by 
$\ov{\nu}_I(t) = \sum_{i \in I} -v(d_i) \etabari$,
if $t$ is induced by the matrix with diagonal entries $d_i$ for all $i \in I$.
Obviously, $\ov{\nu}_I \circ \sigma_I = s_{I} \circ \nu$ on $T$.

Now put
\[ \ov{A} = A \cup \bigcup_{\emptyset \neq I \subset \n} A_I  = 
\bigcup_{\emptyset \neq I \subseteq  \n} A_I.\]
Here of course $A_{\n} = A$ and $s_{\n}$ is the identity.

Let us denote by $\cs_c'(I)$  the set of all seminorms on $V$ which are
canonical with respect to $v_1,\ldots, v_n$, and whose kernel is equal
to $V_{\n \backslash I}$, and by $\cs_c(I)$ the corresponding
quotient space with respect
to equivalence of norms. Then we define a map 
\[\varphi:  A_I \longrightarrow \cs_c(I)\]
by associating to the point $x = \sum_{i \in I} x_i \etabari$ in $A_I$ 
the seminorm represented by 
\begin{eqnarray*}
\ga(\lambda_1 v_1 + \ldots + \lambda_n v_n) =
 \sup \{ |\lambda_i| q^{-x_i}:
i \in I\} 
\end{eqnarray*}
Obviously, $\varphi: A_I \rightarrow \cs_c(I)$ is bijective.

Combining all these maps yields a bijection 
\[\varphi: \ov{A} \longrightarrow \cs_c.\]

Now we want to define an action of $N$ on $\ov{A}$. 
First of all, we use the homomorphism $\ov{\nu}_I \circ \sigma_I = s_I \circ \nu  
: T \rightarrow A_I$ to define an action of $T$ on $A_I$ by
affine transformations. 
For $w \in W$ we denote the induced 
permutation of the set $\n$ also by $w$, i.e. we abuse notation so that 
$w(v_i) = v_{w(i)}$. 
Now define an action of $w$ on $\ov{A}$ by putting together the maps
\[w: A_I \longrightarrow A_{w(I)},\]
sending $\etabari$ to $\ov{\eta}_{w(i)}^{w(I)}$. 
Note that it is compatible with the action of $w$ on $A$, i.e. we have
\[ w \circ s_I = s_{w(I)} \circ w\]
on $A$. 
These two actions give rise to an action of 
$N = T \rtimes W$ on $\abar$, which we denote by $\nu$. 
If $N$ acts in the usual way on $\cs_c$, i.e. by $\gamma \mapsto 
\gamma \circ n\inv$, it can easily be checked that
the bijection $\varphi$ is $N$-equivariant.

Let us now define a topology on $\ov{A}$. For all $I \subset \n$ we put
\[ \Delta_I = \sum_{i \notin I} \R_{\geq 0} \eta_i \subset A.\]
For all open and bounded subsets $U \subset A$ we define
\[ \Gamma_U^I = (U + \Delta_I) \cup \bigcup_{I \subseteq J \subset \n} 
s_J (U + \Delta_I).\]
We take  as a base of our topology on $\ov{A}$ the open subsets 
of $A$ together with these sets  $\Gamma_U^I$ for all non-empty 
$I \subset \n$ and all open bounded subsets $U$ of $A$. 

Note that every point $x \in \ov{A}$ has a countable fundamental 
system of neighbourhoods. This is clear for $x \in A$. 
If $x$ is in $A_I$ for some $I \subset \n$, then choose some 
$z \in A$ with $s_I z = x$, and choose a countable decreasing fundamental 
system of bounded open neighbourhoods $(V_k)_{k \geq 1}$ of $z$ in $A$. 
Put $U_k = V_k + \sum_{i \notin I} k  \eta_i$. This is an open neighbourhood of 
$z +  k \sum_{i \notin I}  \eta_i$. 
Then $(\Gamma_{U_k}^I)_{k \geq 1}$ is a fundamental system of open 
neighbourhoods of $x$.

Our next goal is to compare the previous construction to the compactification
of the appartment $A$ which was used in \cite{we}.
Denote by $V^\ast$ the dual vector space corresponding to $V$. The map 
$g \mapsto ~^t \! g\inv$ induces an isomorphism $\alpha: \underline{PGL}(V^\ast)
\rightarrow \underline{PGL}(V)$, which maps the torus $\bt'$, given 
by the diagonal matrices with respect to the dual basis $v_1^\ast,
\ldots, v_n^\ast$, to the torus 
$\bt$. This defines an isomorphism of the cocharacter groups
\[ \alpha_\ast: X_\ast(\bt') \longrightarrow X_\ast(\bt).\]
Note that $\alpha_\ast$ introduces a sign, i.e. if
$\eta_i'$ is given by the map sending $x$
to the diagonal matrix with $i$-th entry $x$, and entries $1$ at 
the other places, then $\alpha_\ast$ 
maps  $\eta_i'$ to the cocharacter
$-\eta_i$. If we denote by $\Lambda'$ the appartment in the 
building of $\underline{PGL}(V^\ast)$ given by the torus $\bt'$,
the map $\alpha_\ast$ induces an isomorphism of $\R$-vector spaces
$\beta: \Lambda' \rightarrow A$.

Let $I$ be a non-empty subset of $\n$, and let $(V^\ast)_I$ be the subspace
of $V^\ast$ generated by all $v_i^\ast$ for $i \in I$. Then $(V^\ast)_I
 \simeq (V/V_{\n \backslash I})^\ast$. Hence
 $\alpha$ induces an isomorphism 
 $\alpha_I: \underline{PGL}((V^\ast)_I) \rightarrow \underline{PGL}
 (V / V_{\n \backslash I})$ making
 the following diagram commutative

$$
\begin{CD}
\underline{PGL}(V^\ast)^{(V^\ast)_I} @>{\alpha}>> 
\underline{PGL}(V)^{V_{\n \backslash I}}\\
@V{\rho_I}VV @VV{\sigma_I}V\\
\underline{PGL}((V^\ast)_I) @>{\alpha_I}>> \underline{PGL}(V / V_{\n \backslash I})
\end{CD}
$$ 

Restricting $\alpha_I$ to the torus $\bt'_I$ induced by the
diagonal matrices with respect 
to $v_i^\ast$ for $i \in I$, we get an isomorphism
 $\bt'_I \rightarrow \ov{\bt}_I$.
This induces an isomorphism
\[(\alpha_I)_\ast: X_\ast(\bt'_I) \longrightarrow X_\ast(\ov{\bt}_I).\]
If $\eta_i^{'I}$ denotes the cocharacter induced by mapping $x$ to the 
diagonal matrix with entry $x$ at the place $i$, and with entry $1$ at the
places $j \neq i$ in $I$, we have $(\alpha_I)_\ast(\eta_i^{' I}) = 
- \etabari$. Again we put $\Lambda'_I = X_\ast(\bt'_I) \otimes \R$.
Then $(\alpha_I)_\ast$ induces an $\R$-linear isomorphism 
$\beta: \Lambda'_I \rightarrow A_I$. 

Now we put as in \cite{we}, section 3,
 $\Labar' = \Lambda' \cup \bigcup_{\emptyset
\neq I\subset \n} \Lambda'_I$. Then we have defined
 a bijection 
\[ \beta: \Labar' \longrightarrow \abar.\]
Note that by definition $\beta$ is a
homeomorphism,
if $\Labar'$ is endowed with the topology 
of  \cite{we}, section 3. Besides, in \cite{we}, section 3, we defined an action 
$\nu'$ of $N'$ on $\Labar'$, where $N'$ is the normalizer of $\bt'$. Via the isomorphism
$\alpha: N \rightarrow N'$, this is compatible with the $N$-action on $\abar$
defined above. 

Hence we can deduce
\begin{thm}
i) The topological space $\abar$ is compact and contractible, and $A$ is 
an open, dense subset of $\abar$. 

ii) The action $\nu: N \times \abar \rightarrow \abar$ is continuous.
\end{thm}

{\bf Proof: } This follows immediately from \cite{we}, Theorem 3.4 and Lemma
3.5.~

We can endow the space $\cs'_c$ of canonical seminorms with the topology
of pointwise convergence, i.e. with the coarsest topology such that 
for all $v \in V$ the map $\gamma \mapsto \gamma(v)$ from $\cs'_c$ to 
$\R_{\geq 0}$ is continuous. On $\cs_c$ we have the quotient topology.
Both $\cs'_c$ and $\cs_c$ are Hausdorff. 
We will now show that the bijection $\varphi$
defined above is in fact a homeomorphism.

\begin{prop}The $N$-equivariant bijection $\varphi: \abar \rightarrow \cs_c$
is a homeomorphism.
\end{prop}

{\bf Proof: }  Note that it is enough to show that $\varphi$ is continuous, 
since a continuous bijection
from a compact space to a Hausdorff space  is automatically 
a homeomorphism. 

Let $x = \sum_{i \in I} x_i \etabari $ be a point in $A_I$,
and assume that $x_{i_0} = 0$ for some $i_0 \in I$. 
We want to show that $\varphi$ is continuous in $x$. 
If $\gamma$ is the seminorm $\ga(\lambda_1 v_1 + \ldots + \lambda_n v_n) 
= \sup \{ |\lambda_1| q^{-x_1} , 
\ldots, |\lambda_n| q^{-x_n}\}$, then $\varphi(x)$ is represented by $\gamma$.

If $U$ is an open neighbourhood of $\varphi(x) = \{\gamma\}$, 
then it contains a set of the 
form 
\[\{ \{ \beta\} \in \cs_c:  |\beta(v_i) - \gamma(v_i)| <\varepsilon 
\mbox{ for all }i = 1,\ldots, n\}.\]
We find an open and bounded subset $V$ of $A$ so that all $\sum_{i = 1}^n
y_i \eta_i \in V$ satisfy
\begin{eqnarray*} |q^{-y_i+ y_{i_0}}   - q^{-x_i} | < \varepsilon 
& \mbox{ for all } i \in I \mbox{ and}\\
|q^{-y_i + y_{i_0}}|< \varepsilon  & \mbox{ for all } i \notin I.
\end{eqnarray*}
Then $\Gamma_V^I = (V + \Delta_I) \cup \bigcup_{I \subseteq J \subset \n}
s_J(V+ \Delta_I)$
is an open neighbourhood of $x$ in $\abar$. Let us show that $\varphi(\Gamma_V^I)$
is contained in $U$. If $y$ is a point in $\Gamma_V^I$, say
$y = s_J(z)$ for some $I \subseteq J$ and some $z = \sum_{i = 1}^n z_i \eta_i$
in $ V+ \Delta_I$, then 
$\varphi(y)$ is represented by the seminorm $\beta$ with $\beta(\lambda_1
v_1 + \ldots + \lambda_n v_n) 
= \sup \{ |\lambda_i| q^{-z_i + z_{i_0}}: i \in J \}$.
Hence $|\beta(v_i)- \gamma(v_i)| < \varepsilon$ for all $i= 1,\ldots, n$, 
so that $\varphi(y) \in U$. 
Therefore $\varphi$ is indeed continuous.~

\section{Compactification of  the whole building}

Let us first recall the construction of the Bruhat-Tits building
$X$ corresponding to $\bg$.
For every root $a = a_{ij}$ we denote by $U_a$ the root group
in $G$, i.e. the group of 
 matrices $U = (u_{kl})_{k,l}$ such that the 
diagonal elements $u_{kk}$ are equal to one, and all the other entries
apart from  $u_{ij}$ are zero. 
Then we have a homomorphism
\[\psi_a : U_a \longrightarrow \Z \cup \{\infty\}\]
by mapping the matrix $U=(u_{kl})_{k,l}$ to $v(u_{ij})$. Put for all 
$l \in \Z$
\[U_{a,l} = \{ u \in U_a: \psi_a(u) \geq l \}. \]
We also define $U_{a,\infty} = \{1\}$, and $U_{a, -\infty} = U_a$. For all 
$x \in A$ let now $U_x$ be the group generated by 
$U_{a, -a(x)}= \{u \in U_a: \psi_a(u) \geq -a(x)\}$ for all $a \in \Phi$. 
Besides, put $N_x= \{ n \in N: \nu(n)x=x \}$, and 
\[P_x = U_xN_x = N_x U_x.\]
Then the building $X= X(\underline{PGL}(V))$ is given as
\[ X = G \times A /\sim,\]
where the equivalence relation $\sim$ is defined as follows:
\[
(g,x) \sim (h,y),  \begin{array}{l}
\mbox{ iff  there exists an element } n \in N \\
 \mbox{ such that } \nu(n)x= y \mbox{ and } g\inv h n \in P_x.
 \end{array}
 \]
There is a continuous action of $G$ on $X$ via left multiplication on the first 
factor, which  extends the $N$-action on $A$.  
For all $x \in A$ the group $P_x$ is in fact the stabilizer of $x$.

Now we define for all non-empty subsets $\Sigma$ of $\abar$ and all roots 
$a \in \Phi$
\begin{eqnarray*}
 f_\Sigma(a) & =  \inf\{t: \Sigma \subseteq\overline{\{ z \in A: a(z) \geq -t\}}
  \} \\
~ & = -\sup \{ t: \Sigma \subseteq \ov{\{ z \in \La: a(z) \geq t \} } \}
\end{eqnarray*}
Here we put $\inf \emptyset = \sup \R = \infty$ and 
$\inf \R = \sup \emptyset = - \infty$. Obviously, 
$f_x(a) = -a(x)$, if $x$ is contained in $A$.

We define a subgroup 
\[ U_{a,\Sigma} = U_{a,f_\Sigma(a)} = 
\{ u \in U_a : \psi_a(u) \geq f_\Sigma(a)\}\]
of $U_a$, where $U_{a, \infty} = 1$ and $U_{a, -\infty} = U_a$. 
By $U_\Sigma$ we denote the subgroup of $G$ generated by all the $U_{a,\Sigma}$ 
for roots $a \in \Phi$. 

Recall from the previous section
that the isomorphism $\alpha: \underline{PGL}(V^\ast) \rightarrow 
\underline{PGL}(V)$ induces
a homeomorphism $\beta: \Labar' \rightarrow \abar$.
Besides, $\alpha$ induces a bijection between the root systems
$\alpha^\ast:\Phi=    \Phi(\bt) \rightarrow \Phi(\bt')$.
Let us define as in \cite{we}, section 4,
for all 
non-empty subsets $\Omega$ of $\Labar'$ and all roots $a \in \Phi(\bt')$ the subgroup
$U_{a,\Omega}$ of the group $PGL(V^\ast)$ as 
$\{ u \in U_a : \psi_a(u) \geq f_\Omega(a)\}$ for
$f_\Omega(a)  =  
\inf\{t: \Omega \subseteq\overline{\{ z \in \La': a(z) \geq -t\}} \}$.
Besides, let $U_\Omega$ be the subgroup of $PGL(V^\ast)$ generated
by all $U_{a,\Omega}$, and put $P_\Omega = N_\Omega U_\Omega$. 
 Then 
$f_\Omega(\alpha^\ast(a)) = f_{\beta(\Omega)} (a)$, and hence 
$\alpha (U_{\alpha^\ast(a),\Omega}) = U_{a,\beta(\Omega)}$ for all roots $a\in \Phi$, 
which implies
\[\alpha(U_\Omega)= U_{\beta(\Omega)}.\]

Besides, we define for all non-empty $\Sigma \subseteq \abar$ the group
$N_\Sigma = \{n \in N: \nu(n) x = x \mbox{ for all } x \in \Sigma\}$, and 
\[ P_\Sigma= U_\Sigma N_\Sigma = N_\Sigma U_\Sigma.\]
Since $N_\Sigma$ normalizes $U_\Sigma$, which can be shown as in \cite{we}, 4.4,
this set is indeed a group.
Besides, we have
\[\alpha(P_{\Omega}) = P_{\beta(\Omega)}.\]

Now denote by $Z\subset T$ the kernel of the map 
$\nu: T \rightarrow A$. We fix the point $0 \in A$. 
The group $U_0^\wedge = U_0 Z$ is compact and open in $G$, see \cite{la}, 12.12.
 
We define the compactification $\ov{X}$ of the building $X$ as 
\[\ov{X} = U_0^\wedge \times \abar/\sim,\]
where the equivalence relation $\sim$ is defined as follows:
\[
(g,x) \sim (h,y),  \begin{array}{l}
\mbox{ iff  there exists an element } n \in N \\
 \mbox{ such that } \nu(n)x= y \mbox{ and } g\inv h n \in P_x.
 \end{array}
 \]
Since \cite{we}, 4.4 implies that $n P_x n\inv = P_{\nu(n) x}$
for all $n\in N$ and $x \in \abar$, it is easy to check that $\sim$ 
is indeed an equivalence relation.
We equip $\ov{X}$ with the quotient topology. 
Hence $X$ is open and dense in $\ov{X}$. In \cite{we}, 
section 4, we gave a similar definition of a compactification $\ov{X}'$
of 
the building $X'$ corresponding to the group $\underline{PGL}(V^\ast)$,
using the groups $P_x$ for $x \in \Labar'$ recalled above, and the group
$(U_0^\wedge)' = U_0 Z'$ for $Z' = \mbox{ker}(\nu': T' \rightarrow \Lambda')$ and
$0 \in \Lambda'$. Since
 $\alpha(P_x) = P_{\beta (x)}$, the 
homeomorphism 
\[(\alpha,\beta): (U_0^\wedge)' \times \Labar' \longrightarrow 
U_0^\wedge \times \abar\]
induces a homeomorphism
\[ \ov{X}' \longrightarrow \ov{X}.\]
Therefore we can use the results of \cite{we} to deduce
\begin{thm}
$\ov{X}$ is a compact and contractible topological space,
containing the building $X$ as an open, dense subset.
\end{thm}
{\bf Proof: } This follows from \cite{we}, Theorem
5.5 and 5.6.~

For every $x \in \abar$ we have a mixed Bruhat
 decomposition
$G = P_0 N P_x$, which follows from \cite{we}, 4.9. 
This can be used to define an action $\mu$ of $G$ on $\xbar$ in the
 following way: 
 Let $z \in \xbar$ be an element induced by $(u,x) \in U_0^\wedge \times \abar$,
 and let $g \in G$. Write $gu = vnh$ for
 some $v \in U_0^\wedge$, 
 $n \in N$ and $h \in P_x$. Then $\mu (g,z)$ is the point in $\xbar$ induced
 by  $(v, \nu(n) x)$ in $U_0^\wedge \times \abar$.
 Since $n P_x n\inv = P_{\nu(n) x}$,  this induces a well-defined 
 action of $G$ on $\ov{X}$, which extends the action $\nu$ of $N$ on $\abar$. 

Our next goal is to show that the action $\mu$ is continuous. (This is a result
whose counterpart for $\xbar'$ was not proven in \cite{we}.)

\begin{thm} The action $\mu: G \times \xbar \rightarrow \xbar$ defined above is
continuous.
Hence $\xbar$ is equal to the quotient of $G \times \abar$ after the following
equivalence
relation:
\[(g,x) \sim (h,y), \begin{array}{l} \mbox{ iff there exists some }
n \in N \\
\mbox{ with }\nu(n)x=y
\mbox{ and }g\inv h n \in P_x.\end{array} \]
 
\end{thm}

{\bf Proof: } Let us first show that the continuity of the $G$-action $\mu$
implies
that we can also write $\xbar$ as a quotient of $G \times \abar$. 
The inclusion $U_0^\wedge \hookrightarrow G$ induces a 
bijection $\tau: \xbar = ((U_0^\wedge \times \abar) /
 \sim) \rightarrow ((G \times \abar) / \sim)$, which is obviously continuous.
 Besides, if $V$ is an open subset of $\xbar$, then the preimage of $\tau(V)$
 in $G \times \abar$ is equal to $\mu_0\inv(V)$,
where $\mu_0$ is the restriction of $\mu$ to $G \times \abar$.
 Hence the continuity of $\mu$ implies that $\tau(V)$ is open.
 
 If we restrict $\mu$ to the open subgroup $U_0^\wedge$ of $G$, the action
 of $U_0^\wedge$ on $\xbar$ is induced by the natural action 
  of $U_0^\wedge$ on $U_0^\wedge \times \abar$ given by left multiplication
  in the first factor, so that $\mu$ is continuous 
  in all points $(1,z)$ for $z \in \xbar$. 
 
Hence it suffices to show that every $g \in G$ acts continuously on $\xbar$.
Let us denote the quotient map $U_0^\wedge  \times \abar \rightarrow \xbar$
by $\pi$. Then 
we have to check that the map $\mu_g: U_0^\wedge \times \abar 
\rightarrow \xbar$, defined by $\mu_g(u,x) = \pi(v,\nu(n)x)$, if $gu = vnh$ is a
decomposition according to  
$G = U_0^\wedge N P_x$, is continuous.

Assume that $u_k$ is a sequence in $U_0^\wedge$ converging to $u\in U_0^\wedge$,
 and that $x_k$
is a sequence in  $\xbar$ converging to 
$x \in \xbar$. Then we have to show 
that $\mu_{g}(u_k,x_k)$ (or at least a subsequence of it)
converges to $\mu_g(u,x)$. 
Write $g u_k = v_k \gamma_k$ for $v_k \in U_0^\wedge$ and 
$\gamma_k \in N P_{x_k}$. Since $U_0^\wedge$ is compact, we can 
pass to a subsequence of the $(u_k,x_k)$ and assume that 
$v_k$ converges in $U_0^\wedge$ to some element
$v \in U_0^\wedge$. Then the sequence $\gamma_k$  converges in $G$.

Assume for the moment that we have proven the following claim:

{\bf Claim: } If $x_k$ is a sequence in $\abar$ converging to $x \in \abar$, 
and $\gamma_k$ is a  sequence in $N P_{x_k}$ converging to some $\gamma \in G$,
then (after possibly passing to a subsequence of $x_k$) 
we can write $\gamma_k = n_k h_k$ with $n_k \in N$ and $h_k \in P_{x_k}$ such
that $n_k$ converges to some $n \in N$, and $h_k$ converges to some $h \in P_x$.
In particular, $\gamma$ lies in $NP_x$.

Believing this result for a moment we can write 
$\gamma_k = n_k h_k$ with $n_k \rightarrow n \in N$ and $h_k \rightarrow
h \in P_x$. As $g u_k = v_k n_k h_k \in U_0^\wedge N P_{x_k}$ converges
to $gu = v nh \in U_0^\wedge N P_x$, the continuity of the $N$-action
on $\abar$ implies  that  $\mu_g(u_k, x_k) 
= \pi(v_k, \nu(n_k) x_k)$ converges to $\pi(v,\nu(n)x) = \mu_g(u,x)$.

 Hence it remains to prove the claim. As a first step, we will prove it under 
 the condition that the sequence $x_k$ is contained in $A$. The limit point
 $x$ lies in some component $A_I$ for $I \subseteq \n$. We fix an index
 $i_0 \in I$, and write $x_k = \sum_i x_{k,i} \eta_i$
 with $x_{k,i_0} = 0$, and $x = \sum_{i \in I} x_{i} \etabari$ 
 with $x_{i_0} =0$. The convergence
 $x_k \rightarrow x$ translates into 
 \[\begin{array}{ll} x_{k,i} \rightarrow x_i , &
 \mbox{ if } i \in I \mbox{ and}\\
 x_{k,i} \rightarrow \infty, & \mbox{ if } i \notin I.
 \end{array}\]
 
  By the
 pigeon-hole principle there must be one ordering $i_1 \prec i_2 \prec \ldots
 \prec i_r$ of the set $\n \backslash I = \{ i_1, \ldots, i_r\}$ such that 
 infinitely many members of the sequence $x_k$ satisfy 
 \[ x_{k, i_1} \geq x_{k, i_2 } \geq \ldots \geq x_{k, i_r}.\]
 We replace $x_k$ by the subsequence of all $x_k$ satisfying
 these inequalities. 
 Now we choose any ordering $\prec$ of the set $I$ and define a linear
  ordering $\prec$ on the whole
 of $\n$ by combining the orderings on $I$ and $\n \backslash I$ 
subject to the condition $i \prec j$, if $i \notin I$ and $j \in I$. 
Then $\Phi^+ = \{ a_{ij} : i \prec j\}$ is the set of positive roots with 
respect to a suitable chamber corresponding to $\Phi$
(cf. \cite{bou}, chapter V and VI). 

If $a_{ij}$ is contained in  the set $\Phi^-$ of negative roots, 
i.e.  we have $j \prec i$, then 
\begin{eqnarray*}
 f_{x_k}(a_{ij}) = -a_{ij}(x_k) = x_{k,j} - x_{k,i} \left\{
\begin{array}{ll}\geq 0, & \mbox{ if }i,j \notin I\\
\rightarrow \infty, & \mbox{ if } i\in I, j \notin I \\
\rightarrow x_j - x_i, & \mbox{ if } i,j \in I. 
\end{array}\right.
\end{eqnarray*}
Hence there is a real constant $C$ such that 
$f_{x_k}(a) \geq C $ for all $x_k$ and all $a \in \Phi^-$.

Now let us denote by  ${\bf{U}}_{\Phi^+}$, respectively ${\bf{U}}_{\Phi^-}$
 the 
corresponding subgroup of $\bg$ (see \cite{bo}, 21.9), and by
 $U_{\Phi^+}$, respectively $U_{\Phi^-}$, 
 their $K$-rational points. 
For all $y \in \abar$ we put 
\[ U_y^+ = U_{\Phi^+} \cap U_y \quad \mbox{and} \quad U_y^- 
= U_{\Phi^-} \cap U_y.\]
It follows from  \cite{we}, Theorem 4.7, that
the multiplication map induces a bijection 
$\prod_{a \in \Phi^\pm} U_{a,y} \longrightarrow U_y^\pm$, 
where the product on the left hand side may be taken in arbitrary order.
Hence the fact that for all $a \in \Phi^-$ the numbers  $f_{x_k}(a)$ 
are bounded from below implies that all $U_{x_k}^-$ are in fact
contained in a compact subset of $U_{\Phi^-}$. 

It follows from \cite{we}, Corollary 4.8, that the group $P_y$ can be written as 
\[P_y = U_y^- U_y^+ N_y = N_y U_y^+ U_y^-.\]

Therefore the element $\gamma_k  \in N P_{x_k}$ has a product decomposition as
$\gamma_k = n_k u_k^+ u_k^-$ with $n_k \in N$, $u_k^+ \in U_{x_k}^+$, 
and $u_k^- \in U_{x_k}^-$. Since all $u_k^-$ are contained
in a compact subset of $U_{\Phi^-}$, we can pass to a subsequence and
assume that $u_k^-$ 
converges to some element $u^- \in U_{\Phi^-}$. By \cite{we}, Lemma 4.3
and Theorem 4.7, we find that $u^-$ lies in $U_x^-$. 
Besides, after passing to a suitable subsequence of $\gamma_k$, 
all $n_k \in N$ lie in the same coset modulo $T$, i.e. $n_k = m t_k$
for some $m \in N$. Then $t_k u_k^+$ is a converging sequence in the Borel group
$T U_{\Phi^+}$. 
Hence $t_k$ converges to some $t \in T$, and that $u_k^+$ 
converges
to some $u^+ \in U_{\Phi^+}$. Again we deduce from \cite{we},  4.3 and
 4.7, that $u^+$ lies in fact in $U_x^+$. Hence $\gamma_k = m t_k u_k^+
u_k^- = n_k h_k$ for $n_k = m t_k \in N$ and $h_k = u_k^+ u_k^- \in P_{x_k}$,
and we have shown that $n_k$ converges to $mt \in N$ and that $h_k $
converges to $u^+ u^- \in P_x$. Therefore our claim holds if the $x_k$
are contained in $A$. 

In the general case we denote again by $A_I$ the piece of $\abar$ containing
$x$. After passing to a subsequence of $x_k$ we can assume that all $x_k$ lie
in the same piece $A_J$. Then $I$ must be contained in $J$. 
Besides, after passing to a subsequence of $\gamma_k$ we find some $m_0 \in N$ 
such that $m_0 \gamma_k$ is contained in $TP_{x_k}$ for all $k$. 

Now recall the commutative diagram 
$$
\begin{CD}
PGL(V^\ast)^{(V^\ast)_J} @>{\alpha}>> PGL(V)^{V_{\n \backslash J}}\\
@V{\rho_J}VV @VV{\sigma_J}V\\
PGL((V^\ast)_J) @>{\alpha_J}>> PGL(V / V_{\n \backslash J})
\end{CD}
$$ 
from section 3. From the proof of \cite{we}, Theorem 5.7, we know that
for all $y' \in \Lambda'_J$ the  map $\rho_J$ induces a surjection $\rho_J:
P_{y'} \rightarrow P_{y'}^J$, where $P_{y'}^J$ is defined in the same way as 
$P_{y'}$,
replacing $PGL(V^\ast)$ by 
$PGL((V^\ast)_J)$ and the appartment $\Lambda'$ by  $\Lambda'_J$. 
Besides,  $P_y'$ is the full preimage of $P_{y'}^J$ under $\rho_J$.
We know that
for $y = \beta(y')$ the equality 
$\alpha(P_{y'}) = P_{y}$ holds, and it is easily checked that also
$\alpha_J (P_{y'}^J) = \ov{P}_{y}^J$, where for every point $y \in A_J$
the group $\ov{P}^J_y$ is the subgroup of $PGL(V/V_{\n \backslash J})$ 
defined in the same way as $P_y$, replacing the
group $PGL(V)$ by $PGL(V / V_{\n \backslash J})$, 
and the appartment $A$ by $A_J$. 
Hence it follows that $\sigma_J$ induces a surjection $\sigma_J: P_y \rightarrow 
\ov{P}^J_y$, and that $P_y$ is the full preimage of $\overline{P}_y^J$. 

Since $m_0 \gamma_k$ lies in $PGL(V)^{V_{\n \backslash J}}$, we
can apply $\sigma_J$ to this sequence and get a converging sequence
in $PGL(V/V_{\n \backslash J})$. Since $x_k$ is a sequence in $A_J$, we
can apply the case of our claim already proven, this time working
in the building corresponding to $\underline{PGL}(V/V_{\n \backslash J})$.
It follows that $\sigma_J (m_0 \gamma_k)$ can be written as
$\sigma_J (m_0 \gamma_k) = \tilde{m}_k \tilde{h}_k$ with
$\tilde{m_k} \in N(\ov{T}_J)$ converging to some $\tilde{m}$ in 
$N(\ov{T}_J)$ and $\tilde{h}_k \in \ov{P}_{x_k}^J$, converging to 
some $\tilde{h}$ in $\ov{P}_x^J$. 

After passing to a subsequence, we find elements
$m_k$ in $N$ projecting to $\tilde{m}_k$
under $\sigma_J$ such that $m_k$ converges to some $m \in N$ with
$\sigma_J(m) = \tilde{m}$. Then $\sigma_J$ maps $m_k\inv m_0 \gamma_k$ 
to the element $\tilde{h}_k \in \ov{P}^J_{x_k}$. Hence $m_k\inv m_0 
\gamma_k$ lies in $P_{x_k}$. Besides, $m_k\inv m_0 \gamma_k$ converges
to $m\inv m_0 \gamma$, and this projects via $\sigma_J$ to $\tilde{h}$. 
Now $\tilde{h}$ lies in $\ov{P}_x^J$ for the limit point $x \in A_I$, 
and it is easily checked that
$\sigma_J\inv (\ov{P}_x^J)$ is contained in $\sigma_I\inv (\ov{P}_x^I)$,
which is equal to $P_x$. Hence $m\inv m_0 \gamma \in P_x$.
Therefore we can put $n_k = m_0\inv m_k$
and $h_k = m_k\inv m_0 \gamma_k$ to obtain the desired decomposition
of $\gamma_k$: The sequence $n_k$ converges to $n = m_0 m\inv$ in $N$, 
and the sequence $h_k$ converges to  $m\inv m_0 \gamma$ in $P_x$.

This finishes the proof of the theorem.~

\section{$\xbar$ as the space of seminorms on $V$}
We want to extend the homeomorphism $\varphi: \abar \rightarrow \cs_c$, 
where $\cs_c$ is the space of equivalence classes of
canonical seminorms with respect to $v_1,\ldots, v_n$, 
to the 
whole compactified building $\xbar$. Let $\cs'$ denote the set
of all seminorms on $V$, and denote by $\cs$ the quotient of $\cs'$
with respect to the equivalence relation on seminorms. On $\cs'$
we have the topology of pointwise convergence, i.e. the coarsest
topology such that for all $v\in V$ the map $\gamma \mapsto \gamma(v)$
is continuous. The quotient space $\cs$ is equipped with the quotient
topology. Both $\cs'$ and $\cs$ are Hausdorff.
Note that $\cs'$ carries a natural $GL(V)$-action given by
$\gamma \mapsto g(\gamma) = \gamma \circ g\inv$. This induces a 
$G=PGL(V)$-action on the 
quotient space $\cs$. 

\begin{thm}
The map $G \times \abar \rightarrow \cs$, which associates to 
$(g,x) \in G \times
\abar$ the point $g(\varphi(x)) \in \cs$,
induces a $G$-equivariant homeomorphism
\[ \varphi: \xbar \longrightarrow \cs.\]
\end{thm}

{\bf Proof: }In order to show that $\varphi$ is well-defined we have to check 
that for
all $x \in \abar$ the group $P_x$ 
stabilizes the seminorm class $\varphi(x)$. We have $P_x = U_x N_x$.
Since $\varphi$ is $N$-equivariant on $\abar$, the group $N_x$ stabilizes
$\varphi(x)$. It remains to show that for all $a \in \Phi$ the group 
$U_{a,x}$ fixes $\varphi(x)$. Let $A_I$ be the piece of $\abar$ 
containing $x$. Let us first assume that $i \notin I$. 
Every $u \in U_{a_{ij}}$ leaves the vectors $v_l$ for $l \neq j$
invariant and maps $v_j$ to $v_j + \omega v_i$ for some $\omega \in K$.
If $\gamma$ is a seminorm representing $\varphi(x)$, then 
$\gamma \circ u\inv 
(v_j) = \gamma(v_j - \omega v_i) = \gamma(v_j)$, as $\gamma(v_i)=0$. 
Hence every $u \in U_{a_{ij},x}$ stabilizes $\gamma$.
If $ i \in I$ and $j \notin I$, then 
$U_{a_{ij}, x} = 1$, and 
there is nothing to prove. 
It remains to treat the case that both $i$ and $j$ are contained in 
$I$. If $x = \sum_{i \in I} x_i \etabari$, then  every $u \in U_{a_{ij},x}$ 
maps $v_j$
to $v_j + \omega v_i$ for some $\omega \in K$ satisfying $v(\omega) \geq 
f_{x}(a_{ij})=x_j - x_i$. For the usual seminorm  $\gamma$ representing 
$\varphi(x)$ we can 
calculate
\[ \gamma (u\inv (v_j)) = \gamma( v_j - \omega v_i) = \sup\{q^{-x_j},
|\omega| q^{-x_j}\} = q^{- x_j} = \gamma(v_j).\]
As $u$ leaves the other $v_l$ invariant, the group $U_{a_{ij},x}$ fixes 
$\varphi(x)$.

The map $\varphi$ is obviously $G$-equivariant. Since 
 every seminorm on $V$ is canonical with respect
to a suitable basis, $\varphi$ is surjective.

Let us now show that $\varphi$ is injective. 
If $g(\varphi(x)) $ coincides with $h(\varphi(y)) $, it
is easy to see that there exists some $n\in N$ satisfying $ n(\varphi(x) ) 
= \varphi(y)$. Since $\varphi$ is $N$-equivariant and injective
on $\abar$, this implies $\nu(n) x = y$. 
Now $g\inv h n$ is contained in the stabilizer of $\varphi(x)$. 
Hence it remains to show that this stabilizer is equal to  $P_x$. 
If $x$ is contained in $A$, we have the Bruhat decomposition $G = P_x N P_x$.
Hence every $g$ stabilizing $\varphi(x)$ can
 be written as $g = p n q$ with $p$ and $q$ in 
$P_x$. We already know that $P_x$ is contained in the stabilizer of $\varphi(x)$, 
hence it follows that $n$ 
stabilizes $\varphi(x)$. Since $\varphi: \abar \rightarrow \cs_c$ is 
$N$-equivariant and injective, 
this implies that $n$, and hence also $g$, is contained in $P_x$.
If $x$ is a boundary point, i.e. $x \in A_I$ for some $I \subset \n$, then 
$\varphi(x)$ induces an equivalence class of norms on the quotient space 
$V / V_{\n \backslash I}$. Every element $g$ stabilizing $\varphi(x)$ leaves
$V_{\n \backslash I}$ invariant, 
so that we can apply $\sigma_I: PGL(V)^{V_{\n \backslash
I}} \rightarrow PGL(V/V_{\n \backslash I})$ to $g$. By the first case, applied to 
the appartment $A_I$ in the building associated to $\underline{PGL}
(V/V_{\n \backslash I})$, 
we find that $\sigma_I(g)$ lies in the group $\ov{P}_x^I$, i.e. the subgroup 
of $PGL(V/V_{\n \backslash I})$ defined in the same way as $P_x$, replacing
$PGL(V)$ by $PGL(V/V_{\n \backslash I})$ and $A$ by $A_I$. We have seen in 
the proof of 4.2 that $\sigma_I\inv (\overline{P}_x^I) = P_x$, so that 
$g$ lies indeed in $P_x$.

Let us now show that $\varphi$ is continuous. 
Since $\varphi$ is induced by the composition 
\[ G \times \abar \stackrel{\mbox{\small id} \times \varphi}{\longrightarrow}
G \times \cs_c \longrightarrow \cs,\]
we have  to show that the  second map
is continuous. 
Hence we have to check that the  action $GL(V) \times \cs'_c \rightarrow \cs'$
is continuous, which amounts to checking that for all $v\in V$ the map
$\psi_v: GL(V) \times \cs'_c \rightarrow \R$, given by 
$(g,\gamma) \mapsto \gamma(g\inv v)$
is continuous. 

We claim that for all $v \in V$, $\gamma_0 \in \cs'_c$ and $\varepsilon >0$ 
there exists an open neighbourhood 
$W$ of $v$ in $V$ and an open neighbourhood $\Gamma$ of $\gamma_0$ in $\cs_c'$
such that for all $w\in W$ and all $\gamma \in \Gamma$ we have
$|\gamma(w) - \gamma_0(w)| < \varepsilon$. 

Let us believe the claim for a second.
Then we find for  $v$ in $V$  and $\gamma_0$ 
in $\cs_c'$ and for every $\varepsilon>0 $ open 
neighbourhoods
$H$ of $1$ in $GL(V)$ and $\Gamma$ of $\gamma_0$ in $\cs_c'$ 
such that for all $h \in H$ 
 and  all $\gamma \in \Gamma$ the estimate $|\gamma(h\inv v) - \gamma_0
 (h\inv v)| < \varepsilon$ holds. Besides, we can make $H$ so small 
that all $h \in H$ satisfy
 $|\gamma_0(h\inv v) - \gamma_0(v)| < \varepsilon$.
Hence
\[|\gamma(h\inv v) - \gamma_0(v)| \leq |\gamma(h\inv v) - \gamma_0(h\inv v) |
+ | \gamma_0 (h\inv v) - \gamma_0(v)| < 2 \varepsilon.\]
This shows that  $\psi_v$ is continuous in $(1,\gamma_0)$, hence everywhere.

It remains to show the claim. 
Let $v = \sum_{i = 1}^n \mu_i v_i$.
Let us first assume that $\gamma_0(v) \neq 0$. Choose some $0 <\delta < 
\min\{1, |\mu_i|: \mu_i \neq 0\}$ such 
that $\delta \gamma_0(v_i) \leq \gamma_0(v)$ for all $i = 1,\ldots, n$.
Then we put 
$W = \{ w = \sum_i \lambda_i v_i: |\lambda_i - \mu_i| < \delta\}$.
For all $w = \sum_i \lambda_i v_i \in W$ 
the conditions on $\delta $ imply that $|\lambda_i| = |\mu_i|$ for 
all $i$ such that $\mu_i \neq 0$. If on the other hand $\mu_i =0$, then 
$|\lambda_i| < \delta$, which implies $|\lambda_i| \gamma_0(v_i) \leq \gamma_0 (v)$.
Hence for all $w \in W$
\begin{eqnarray*} \gamma_0(w) &= 
&  \sup\{|\lambda_1| \gamma_0(v_1), \ldots, |\lambda_n| \gamma_0(v_n)\}\\
 ~& =&  \gamma_0(v)
 \end{eqnarray*}
and
\[ \gamma(w) = \sup \{ \gamma(v), \sup \{ |\lambda_i| \gamma(v_i): \mu_i =0\}\}.\]
Let us denote by $\Gamma$ the open neighbourhood of $\gamma_0$ in $\cs_c'$ 
consisting
of all $\gamma$ satisfying $|\gamma(v)- \gamma_0(v)| < \varepsilon$ and 
$|\gamma(v_i)
- \gamma_0(v_i)| < \varepsilon$ for $i = 1,\ldots,n$.

If $\gamma \in \Gamma$ satisfies
$\gamma(w) = |\lambda_i| \gamma(v_i)|$ for some $i$ with $\mu_i =0$, then 
\begin{eqnarray*}
 \gamma(v) & \leq &|\lambda_i| \gamma(v_i) 
   \leq |\lambda_i | (\gamma_0(v_i) + \varepsilon)
   \leq \gamma_0(v) + \varepsilon \\
 & < & \gamma(v) + 2 \varepsilon,
 \end{eqnarray*}
 which  implies
\begin{eqnarray*} |\gamma(w) - \gamma_0(w)| & = &
|\, |\lambda_i| \gamma(v_i) - \gamma_0(v) |\\ 
&  \leq & |\, |\lambda_i| \gamma(v_i) - \gamma(v)| + 
| \gamma(v) - \gamma_0(v)|\\
 & < & 3 \varepsilon.
 \end{eqnarray*}
If $\gamma(w) = \gamma(v)$, this estimate holds trivially, so that 
our claim follows. 
 
Now we treat the case that $\gamma_0(v) = 0$, i.e. we have $\gamma_0(v_i) = 0$ for 
all $i$ such that $\mu_i \neq 0$. 
Choose some $0 <\delta < 
\min\{1, |\mu_i|: \mu_i \neq 0\}$, put
 $W = \{ w = \sum \lambda_i v_i: |\lambda_i - \mu_i| < \delta\}$, and let
$\Gamma$ be the open neighbourhood of $\gamma_0$ consisting of all $\gamma$ satisfying
$|\gamma(v)| < \varepsilon$ and
$|\gamma(v_i) - \gamma_0(v_i)| < \varepsilon $ for all $i = 1,\ldots, n$.
As above, we have for all $w \in W$ 
\[\gamma(w) = \sup\{ \gamma(v), \sup\{|\lambda_i| \gamma(v_i): \mu_i = 0\}\}.\]
If $\gamma \in \Gamma$ satisfies $\gamma(w) = \gamma(v)$, then
for all $i$ such that $\mu_i = 0$ the estimate $|\lambda_i| \gamma(v_i) 
\leq \gamma(v) <
\varepsilon$ holds. Since $\gamma_0(w) = |\lambda_i| \gamma_0(v_i)$ for some
$i$ satisfying $\mu_i = 0$, this implies that 
\begin{eqnarray*} |\gamma(w) - \gamma_0(w)| & = &
|\gamma(v) - |\lambda_i| \gamma_0(v_i)|\\
& \leq & |\gamma(v)- |\lambda_i| \gamma(v_i)| + |\lambda_i||\gamma(v_i) -\gamma_0(v_i) |
< 2 \varepsilon.
\end{eqnarray*}.

If $\gamma \in \Gamma$ satisfies $\gamma(w) = |\lambda_i| \gamma(v_i)$ 
for some $i$ with  
$\mu_i = 0$, then we have for all $j$ such that $\mu_j = 0$:
\begin{eqnarray*} |\lambda_j| \gamma_0(v_j) &\leq& |\lambda_j| \gamma(v_j) + 
\varepsilon \\
&\leq& |\lambda_i| \gamma(v_i) + \varepsilon \leq |\lambda_i| \gamma_0(v_i)
+ 2 \varepsilon. 
\end{eqnarray*}
If $\gamma_0(w) = |\lambda_j| \gamma_0(v_j)$, we also have 
$|\lambda_i| \gamma_0(v_i) \leq
|\lambda_j| \gamma_0(v_j)  $, whence
\[|\gamma(w) - \gamma_0(w)| = ||\lambda_i| \gamma(v_i) - |\lambda_j| \gamma_0(v_j)|
\leq \varepsilon .\]
Hence the claim is proven.

Therefore $\varphi$ is a continuous bijection from a compact space to 
a Hausdorff space, hence it is a homeomorphism. ~

\section{The reduction map from Drinfeld's symmetric domain}
Let $\PP(V) = \mbox{Proj Sym}\, V$ be the projective space corresponding
to $V$, i.e. points in $\PP(V)$ correspond to lines in the dual space of $V$. 
Drinfeld's $p$-adic symmetric domain $\Omega$ is the  complement in $\PP(V)$ 
of the union of all $K$-rational hyperplanes. $\Omega$ carries the structure
of a rigid analytic variety. There is a reduction map 
$r: \Omega \rightarrow X$ from $\Omega$ onto the building $X$, see \cite{dr},
\S 6, 
which is defined as follows: 
Every $\ov{K}$-rational point $x$ in $\Omega$ induces a line in the 
 dual space of $(V \otimes \ov{K})$. For every element $z \neq 0$ on this line, 
 the map 
$v \mapsto |z(v)|_{\overline{K}}$ defines a norm on $V$. 
Then $r(x)$ is the point in $X$
associated  to this norm via the bijection $\varphi$ discussed in the 
previous section. 

We want to extend this reduction map to a map from the whole projective space
to our compactification $\xbar$ of the building. In the following, 
we identify $\xbar$ with the set of equivalence classes of 
seminorms on $V$ without specifying the homeomorphism $\varphi$ any longer.

Also we consider Berkovich spaces instead of rigid analytic varieties. 
Namely, let $\PP(V)^{an}$ and $\Omega\an$ be the analytic spaces in the sense
of \cite{be1}
corresponding to the projective space and the $p$-adic symmetric
domain. Then $\PP(V)\an$ can be identified with the set of equivalence
classes of multiplicative seminorms on the polynomial ring $\mbox{Sym}\, V$
extending the absolute value on $K$, which do not vanish identically on $V$. 
Here two such seminorms $\alpha$ and $\beta$
are equivalent, iff there exists a constant $c> 0$ such that for all 
homogeneous polynomials $f$ of degree $d$ we have $\alpha(f) = c^d \beta(f)$,
see \cite{be2}.  

In \cite{be2}, Berkovich defines a continuous, $PGL(V)$-equivariant reduction
map $r : \Omega\an  \rightarrow X$, and a right-inverse
$j: X \rightarrow \Omega\an$, which identifies $X$ homeomorphically with a closed
subset of $\Omega\an$. 

We can now generalize these results to the compactifications $\PP(V)\an$ and $\xbar$
by using almost verbatim the same constructions for $r$ and $j$ as in \cite{be2}. 

For every point in $\PP(V)\an$ represented by the seminorm $\alpha$ on the
polynomial ring $\mbox{Sym} \, V$, the restriction of $\alpha$ to $V$ induces
a seminorm on $V$, hence a point in $\xbar$. This induces a 
$PGL(V)$-equivariant map 
\[r: \PP(V)\an \rightarrow \xbar.\]

On the other hand, let $x$ be a point in $\xbar$, corresponding to the class of 
seminorms on $V$ represented by $\gamma$. Then $\gamma$ is canonical with respect to
a basis  $w_1, \ldots, w_n$  of $V$. If $f = \sum_{\nu = (\nu_1,\ldots, \nu_n)}
a_\nu w_1^{\nu_1} \ldots w_n^{\nu_n}$ is a polynomial in $\mbox{Sym} \, V$, we put
\[
\alpha(f) = \sup\{|a_\nu| \gamma(w_1)^{\nu_1} \ldots  \gamma(w_n)^{\nu_n} \}.\]
Then $\alpha$ is a multiplicative seminorm on $\mbox{Sym} \, V$ extending
the absolute value of $K$, hence it induces a point in $\PP(V)\an$.

This defines a $PGL(V)$-equivariant map
\[ j : \xbar \longrightarrow \PP(V)\an.\]

\begin{prop}
The maps $r: \PP(V)\an \rightarrow \xbar$ and $j: \xbar \rightarrow \PP(V)\an$
are continuous and satisfy $r \circ j = \mbox{id}_{\xbar}$.
Besides, $j$ is a homeomorphism from $\xbar$ to its image  $j(\xbar)$,
which is a closed subset of $\PP(V)\an$.
\end{prop}

{\bf Proof: }The continuity of $r$ follows immediately from the definitions.
Since $j$ is obviously continuous on $\abar$ and $PGL(V)$-equivariant,
it is continuous on the whole of $\xbar$. 
By construction, we have $r \circ j = \mbox{id}_{\xbar}$,
so that $j$ is a homeomorphism onto its image, which is a closed subset
of $\PP(V)\an$, since $\xbar$ is compact and $\PP(V)\an$ is Hausdorff.~

\end{document}